\documentclass[11pt]{amsart}
\usepackage{amscd}
\usepackage[all]{xy}
\usepackage{graphicx}
\usepackage{amsmath}
\usepackage{color}

\usepackage{geometry}                
\usepackage[parfill]{parskip}    
\usepackage{amsmath, latexsym, amssymb}
\usepackage{epstopdf}
\numberwithin{equation}{section}
\theoremstyle{plain}

\theoremstyle{definition}

\begin{document}
\newcommand{\R}{{\mathbb R}}
\newcommand{\C}{{\mathbb C}}
\newcommand{\F}{{\mathbb F}}
\renewcommand{\O}{{\mathbb O}}
\newcommand{\Z}{{\mathbb Z}} 
\newcommand{\N}{{\mathbb N}}
\newcommand{\Q}{{\mathbb Q}}
\renewcommand{\H}{{\mathbb H}}

\newcommand{\Aa}{{\mathcal A}}
\newcommand{\Bb}{{\mathcal B}}
\newcommand{\Cc}{{\mathcal C}}    
\newcommand{\Dd}{{\mathcal D}}
\newcommand{\Ee}{{\mathcal E}}
\newcommand{\Ff}{{\mathcal F}}
\newcommand{\Gg}{{\mathcal G}}    
\newcommand{\Hh}{{\mathcal H}}
\newcommand{\Kk}{{\mathcal K}}
\newcommand{\Ii}{{\mathcal I}}
\newcommand{\Jj}{{\mathcal J}}
\newcommand{\Ll}{{\mathcal L}}    
\newcommand{\Mm}{{\mathcal M}}    
\newcommand{\Nn}{{\mathcal N}}
\newcommand{\Oo}{{\mathcal O}}
\newcommand{\Pp}{{\mathcal P}}
\newcommand{\Qq}{{\mathcal Q}}
\newcommand{\Rr}{{\mathcal R}}
\newcommand{\Ss}{{\mathcal S}}
\newcommand{\Tt}{{\mathcal T}}
\newcommand{\Uu}{{\mathcal U}}
\newcommand{\Vv}{{\mathcal V}}
\newcommand{\Ww}{{\mathcal W}}
\newcommand{\Xx}{{\mathcal X}}
\newcommand{\Yy}{{\mathcal Y}}
\newcommand{\Zz}{{\mathcal Z}}

\newcommand{\zt}{{\tilde z}}
\newcommand{\xt}{{\tilde x}}
\newcommand{\Ht}{\widetilde{H}}
\newcommand{\ut}{{\tilde u}}
\newcommand{\Mt}{{\widetilde M}}
\newcommand{\Llt}{{\widetilde{\mathcal L}}}
\newcommand{\yt}{{\tilde y}}
\newcommand{\vt}{{\tilde v}}
\newcommand{\Ppt}{{\widetilde{\mathcal P}}}
\newcommand{\bp }{{\bar \partial}} 

\newcommand{\Remark}{{\it Remark}}
\newcommand{\Proof}{{\it Proof}}
\newcommand{\ad}{{\rm ad}}
\newcommand{\Om}{{\Omega}}
\newcommand{\om}{{\omega}}
\newcommand{\eps}{{\varepsilon}}
\newcommand{\Di}{{\rm Diff}}
\newcommand{\im}{{\rm Im}}
\newcommand {\sppt}{{\rm sppt}}
\newcommand{\Pro}[1]{\noindent {\bf Proposition #1}}
\newcommand{\Thm}[1]{\noindent {\bf Theorem #1}}
\newcommand{\Lem}[1]{\noindent {\bf Lemma #1 }}
\newcommand{\An}[1]{\noindent {\bf Anmerkung #1}}
\newcommand{\Kor}[1]{\noindent {\bf Korollar #1}}
\newcommand{\Satz}[1]{\noindent {\bf Satz #1}}

\renewcommand{\a}{{\mathfrak a}}
\renewcommand{\b}{{\mathfrak b}}
\newcommand{\e}{{\mathfrak e}}
\renewcommand{\k}{{\mathfrak k}}
\newcommand{\pg}{{\mathfrak p}}
\newcommand{\g}{{\mathfrak g}}
\newcommand{\gl}{{\mathfrak gl}}
\newcommand{\h}{{\mathfrak h}}
\renewcommand{\l}{{\mathfrak l}}
\newcommand{\sm}{{\mathfrak m}}
\newcommand{\n}{{\mathfrak n}}
\newcommand{\s}{{\mathfrak s}}
\renewcommand{\o}{{\mathfrak o}}
\newcommand{\so}{{\mathfrak so}}
\renewcommand{\u}{{\mathfrak u}}
\newcommand{\su}{{\mathfrak su}}
\newcommand{\ssl}{{\mathfrak sl}}
\newcommand{\ssp}{{\mathfrak sp}}
\renewcommand{\t}{{\mathfrak t }}
\newcommand{\Cinf}{C^{\infty}}
\newcommand{\la}{\langle}
\newcommand{\ra}{\rangle}
\newcommand{\half}{\scriptstyle\frac{1}{2}}
\newcommand{\p}{{\partial}}
\newcommand{\notsub}{\not\subset}
\newcommand{\iI}{{I}}               
\newcommand{\bI}{{\partial I}}      
\newcommand{\LRA}{\Longrightarrow}
\newcommand{\LLA}{\Longleftarrow}
\newcommand{\lra}{\longrightarrow}
\newcommand{\LLR}{\Longleftrightarrow}
\newcommand{\lla}{\longleftarrow}
\newcommand{\INTO}{\hookrightarrow}

\newcommand{\QED}{\hfill$\Box$\medskip}
\newcommand{\UuU}{\Upsilon _{\delta}(H_0) \times \Uu _{\delta} (J_0)}
\newcommand{\bm}{\boldmath}




\title[Compact symplectic manifolds
of low cohomogeneity]{\large  Compact symplectic manifolds
of low cohomogeneity (corrected version) }
\author[H. V. L\^e  ] {H\^ong V\^an  L\^e
  $^*$}
\address{Institute  of Mathematics of ASCR,
Zitna 25, 11567  Praha 1, Czech Republic} 
\thanks{partially supported by RVO: 67985840} 
\maketitle




\abstract   This  is a   corrected  version of  my  paper published  in  Journal  of Geometry and Physics  25(1998), 205-226.   I added    missing cases   to the  classification theorem 1.1, namely  the   $SO(n+1)$-manifold  $SO(n +2)/ (SO (n) \times SO(2))$, 
the $SO(3)$-manifold $\C P^2$  and  the $SU(3)$-manifold 
$\C P^1 \times \C P^1$.
\endabstract

{\bf Preface}
 Christopher T. Woodward, in his review 
MR1619843  in MathSciNet,  pointed out  a  gap  in  the classification   of  compact symplectic   manifolds   of  cohomogeneity one in my paper   \cite{Le1998}.     ``Unfortunately, there is a mistake in (2.4). The author assumes that the map $\mu:M\to \triangle$, where $\triangle$ is the moment polytope, is smooth. This is not the case, for example, for M the product of two projective lines, and G=SU(2) acting diagonally. Therefore, his conclusions are only valid under this assumption.''  The aim of this   version is to correct  that  mistake  and    to find the missing cases in the  previous classification. I also   slightly improved   the exposition of the previous      version  by adding  three footnotes, inserting    few explanations,   four   new references (including   the    previous version of this paper),     deleting some  unimportant  and imprecise remarks in the previous   version and  polishing   few   sentences.  The  main correction   concerns the classification theorem  1.1.  I have added   Corollary 1, Lemma 2, 
Proposition 4, Lemma 5 and   relations (E1), (E2), (E3), (E4), (E5), (E6) in the new version and  modified   the previous  Proposition 2.3. In the  revised version  the   previous formula (2.4)   and Proposition 2.3 (now is Lemma 3)   are  applied only to  special cases.

\section{ Introduction }

An action of a Lie group $G$  on a manifold $M$ is
 called of {\it cohomogeneity $k$} if
the regular (principal) $G$-orbits  have codimension $k$
in $M$. In other words the orbit space $M/G$ has dimension
$k$.
It is well-known (see e.g. \cite{Kir}) that homogeneous symplectic manifolds are  locally
symplectomorphic to
 coadjoint orbits of Lie groups  whose
symplectic geometry
can be investigated in many aspects  \cite{Gr, H-V, G-K}.
Our motivation is to find a wider class of symplectic manifolds
via group approach, so that they could serve as test examples
for many questions in symplectic geometry  (and symplectic topology).  In this
 note we describe all compact symplectic manifolds admitting a
Hamiltonian action with cohomogeneity 1 of a compact Lie group. We always
assume that the
action is effective. We also
remark that  4-manifolds admitting
symplectic group
actions (of cohomogeneity 1 or
of $S^1$-action) have been studied intensively by many authors, see \cite{Au}
for references. In particular the classification
of compact symplectic 4-manifolds admitting $SO(3)$-action
of cohomogeneity 1 was done by Iglesias \cite{I}.
\medskip

Let us recall that if an action of a Lie group $G$  on $(M,\om)$
preserves the symplectic form $\om$ then there
is a  Lie algebra homomorphism
$$ g = Lie\, G \stackrel{\Ff_*}{\to} Vect_{\om}(M), \eqno (1.1)$$
where $Vect_{\om}(M)$ denotes  the Lie algebra of symplectic vector fields.  The action of $G$ is said to be {\it almost Hamiltonian}
if the image of $\Ff_*$ lies in the subalgebra $Vect_{Ham}(M)$
of Hamiltonian vector fields.  Finally, if the
map $\Ff_*$ can be lifted to a homomorphism $ g \stackrel{\Ff}
{\to}C^{\infty}(M,\R)$ (i.e. $\Ff _*v = sgrad\, \Ff_v$) then the action of $G$ is
called {\it Hamiltonian}. 
  In this note we shall prove the following theorem.

\medskip

{\bf Theorem 1.1}. (Corrected) \footnote{The    original classification  \cite[Theorem 1.1]{Le1998}  corresponds to  the case (i)  in  Theorem 1.} {\it Suppose that a compact symplectic manifold $(M, \om)$ is provided with a Hamiltonian action
  of a compact Lie group $G$ such that  dim $M/G$ =1. Let $\mu$ denote  the   moment map of the $G$-action. \\
1)  If  $\dim  \mu ^{-1}  (m) \le  1$ for all $m \in  M$ then $M$ is  $G$-diffeomorphic either to a $G$-invariant bundle
over a coadjoint
orbit of $G$ whose fiber is a complex projective manifold, or
to a symplectic blow-down of such a bundle along two singular $G$-orbits.

2)  If  there  is a point $m$ such that $\dim \mu^{-1} (m) \ge 2$, then  $M$ is  a  direct  product of a coadjoint orbit of $G$ with one of the following   symplectic    $G$-manifolds of cohomogeneity 1:    the  Grassmannian $SO(n+2)/ (SO(n) \times SO(2))$  with the canonical  Hamiltonian  $SO(n+1)$-action, the  $SO(3)$-manifold $\C P^2$  with  the Hamiltonian action  of $SO(3)$ via the embedding $SO(3) \to SU(3)$, or  the product   $\C P^1 \times \C P^1$  with  the diagonal Hamiltonian action of $SU(3)$.}

\medskip

 The main ingredients of the proof of Theorem 1.1
are the existence of the moment map, {\it the classification of  coadjoint orbits  of  compact Lie group (Table A.3 in Appendix A),    the  classification  of  Riemannian    manifolds of cohomogeneity 1 due  to Alekseevskii-Alekseevskii \cite{AA1993}}, Duistermaat's-Heckman's theorem \cite{D-H}, the convexity theorem of Kirwan \cite{Kiw}. For
certain $G$-diffeomorphism types of these spaces we shall give a complete classification up to
equivariant symplectomorphism (see Section 2).

\medskip

In section 3 we give
a computation of the (small) quantum cohomology ring
of some  spaces admitting a Hamiltonian $U_n$-action
with cohomogeneity 1 and discuss its corollaries.

\medskip

We also consider the case of a symplectic
action of cohomogeneity 2. In particular we get:

\medskip

{\bf Theorem 1.2}. {\it Suppose that a compact symplectic manifold $M$ is provided with a Hamiltonian action
of a compact Lie group $G$ such that dim $M/G$=2.
Then all the principal orbits of $G$ must be 
either (simultaneously) coisotropic or (simultaneously)
symplectic. Thus a principal
orbit of $G$ is either diffeomorphic
to a $T^2$-bundle over a  coadjoint
orbit of $G$ (in the first case) or diffeomorphic to a coadjoint
orbit of $G$ (in the second case).}

\medskip

At the end of our
note we collect in Appendix  A some useful facts of
the symplectic structures on the
coadjoint orbits of compact Lie groups.


\section {Classification of compact symplectic manifolds admitting a
Hamiltonian action with cohomogeneity 1 of a compact Lie group.}


It is known [Br] that if an action of a compact Lie group $G$ on a compact
oriented manifold $M$ has cohomogeneity 1  (i.e. dim $M/G$ =1) then the topological
space $Q= M/G= \pi (M)$ must be either diffeomorphic to the interval $[0,1]$ or a circle $S^1$.
The slice theorem gives us immediately that $G(m)$ is
a principal orbit if and only if the image $\pi (G(m))$
in $Q$ is a interior point. In what follows we assume that $(M,\om)$ is
symplectic and the action of $G$ on $M$ is Hamiltonian.
Under this assumption the quotient $Q$ is  [0,1] (Proposition  2.3).
\medskip

To study  $G$-action on $(M, \om)$   it is  useful   to   fix  a $G$-invariant compatible  metric on $M$, whose existence is well-known,  see e.g.  \cite[Proposition 2.50]{McD-S}.

 {\bf Proposition 2.1.} {\it Let $G(m)$ be a principal
orbit of a Hamiltonian $G$-action on $(M^{2n},\om)$. Then $G(m)$ is
a $S^1$-bundle over a  coadjoint orbit of $G$. }

\medskip

{\it Proof}. 
Let us consider the moment map 

$$ M^{2n} \stackrel{\mu}{\to }g^*:\; \la \mu(m), w\ra = \Ff_w (m).\eqno (2.1)$$
For a vector $V \in T_*G(m)$ there is a vector $v\in g$ such that
$V  = {d\over dt}_{t=0} (exp \, tv)= sgrad\Ff_v$. Hence we get
$$\la \mu_*(V),w\ra = d\Ff_w (V) =\{ \Ff_w, \Ff_v \}(m)= \la [w,v], \mu(m)\ra \eqno (2.2) $$
which implies that  $\mu$ is an equivariant map. Therefore
the image $\mu (G(m))$ of any orbit $G(m)$ on $M$ is a   coadjoint orbit $G(\mu(m))\subset g^*$.

To complete the proof of Proposition 2.1 we look at the
preimage $\mu^{-1}\{ \mu (m)\}$.

\medskip

{\bf Lemma 2.2.} {\it Let $G(m)$ be  a principal orbit. Then the preimage $\mu^{-1}\{ \mu(m)\}$
is  a (connected) orbit of a connected subgroup $S^1_m\subset G$.}

\medskip

{\it Proof.} Clearly the
preimage is a closed subset. 
Let $V$ be a non-zero tangent vector to
the preimage $ \mu ^{-1} \{ \mu (m)\}$ at $x$. Then $ \mu_*(V)=0$. Using the formula
$$ \la \mu_* (V), w\ra = d\Ff_w(V) = \om (sgrad \, \Ff_w, V)\eqno (2.3)$$
for  all $w  \in g$ 
we conclude that $V$ is also a tangent vector to $G(m)$ and moreover
$V$ annihilates the
space $span \, \{ d\Ff _a \, | a \in g\}$ which has codimension 1 in $T^*M$.
In particular there
is an element $\bar{v}\in g$ such that $ V= sgrad\, \Ff_{\bar v}$.
Our claim  on the manifold structure now follows
from the fact that
$exp \, t\bar{v} (x) \subset \mu ^{-1} \{ \mu (m)\}$. 
Finally, the preimage 
is connected because the quotient $\mu (G(m)) = G(m) / \{ \mu^{-1} (\mu(m))\}$
is simply-connected (see Appendix A) and $G(m)$ is connected.\QED

Clearly Lemma 2.2  yields
Proposition 2.1.\QED



We obtain immediately from  Proposition 2  the following.

{\bf Corollary  1.} {\it For any $m \in M$  in a principal $G$-orbit  there exists  a $G$-invariant symplectic  2-form $\bar \om$ on the coadjoint orbit $G (\mu (m))$ such that  $\om|_{G(m)} =  \mu ^* (\bar \om)$.}


{\bf Proposition  2.3.} {\it  The quotient space $M/G$ is   $[0,1]$.}

\begin{proof}  Assume the opposite, i.e.   $M/G = S^1$.  In this case  it is well-known that  the projection
$M \to M/G$  is a   fibration whose   fiber  is  a principal    orbit $G(m ) = G/ G_m$    and whose  structure  group    is  $N_G (G_m)$  \cite[Theorem 5.8]{Br}.  Here we denote by $G_m$ the stabilizer
of  $m \in M$   and by $N_G (G_m)$  the normaliser
of $G_{m}$ in $G$.  Hence,   given a point $m_0 \in M$,     we have the following $G$-equivariant identification \cite{A-A, Br}
$$ M = \R \times _h G_{m_0},$$
where $(t, gG_{m_0})$ is identified with $(t + 1, g h G_{m_0})$ for some element $h\in N_G(G_{m_0})$.

 Denote by $J$  the    $G$-invariant    almost complex structure     that is associated
with   the  given $G$-invariant  compatible  metric on $M$. For each $m \in   M$ denote by
$V_m$  the unit   vector that tangent   to  $\mu^{-1}(\mu(m))$ at $m$.  Since  $\om (V_m, T_m  G(m)) = 0$, it follows that
$\la JV_m,    T_m G(m) \ra = 0$.  Now we choose  the orientation of $V_m$    such that  $\pi_*(JV_m) = \p t$. 
Corollary   1    implies that  the   symplectic   form $\om$ on $M = R\times _h G(m_0)$ has the following form
$$\om (t, y)  = \mu ^*(\bar \om  (t, y))  +  g(t) dt \wedge \alpha \eqno( E 1), $$ 
where $\bar \om  (t, y)$  is  a $G$-invariant   symplectic  form on $\mu (G(m_0))$,  $\alpha$ is the    $G$-invariant    connection 1-form on   $G(m_0)$  associated
 with the  principal  $S^1$-bundle  $G(m_0) \to  \mu(G(m_0))$  and $g(t) \not = 0$, where  the $S^1$  action at $m$ is generated by $\exp V_m$.  

The closedness  of  $\om$ implies that  
$$ \p _t  \bar  \om(t, y)  =   g(t) \cdot d \alpha . \eqno ( E  2)$$
Since  the cohomology class  $[d\alpha ]$  is  the Chern class  of the  $S^1$-bundle $G(m_0) \to \mu  (G(m_0))$,  which  does not depend on $t$, the  equality 
(E2) holds      only if  $[d\alpha]  = 0$, i.e.   the $U(1)$-bundle $G(m_0)  \to  \mu (G(m_0))$ is the trivial   bundle. Hence, denoting $v =   \mu (m_0)$,    we  have
$$ M =   T^2 \times G/ Z(v), \eqno ( E3) $$
where $Z(v)$ is the stabilizer of $v$. 
The closedness of $\om$ implies $\p_t  \bar  \om(t, y) = 0$. Hence  we obtain the following

{\bf Lemma 2.} {\it  $(M, \om)$ is $G$-symplectomorphic  to $(G/Z(v)\times  T^2, \bar \om + dt \wedge  \alpha)$  where $\bar \om$ is a $G$-invariant symplectic form  and $\alpha$  is the canonical  1-form
on the second  fact  $S^1$ of $T^2$.}

 Let $R:  C^\infty (M) \to  C^\infty (T^2)$  be defined  as follows
$$ R(h) (t_1, t_2) : =  h(m_0, t_1, t_2).$$
By Lemma 2, it is not hard to see that  the  composition $R \circ \Ff : g \to C^\infty (T^2)$  defines   a Hamiltonian action
 of $G$ on  $T^2$, which is   of cohomogeneity 1.     It is well-known that  there is no such an  Hamiltonian action
 on torus $T^2$, see  e.g. \cite{Au}.
 Hence  follows Proposition 2.3. 
\end{proof}




{\bf Remark 2.5.}  Let   $m$ belong to a principal  orbit  and $G_m$   its stabilizer.  Then the stabilizer $Z(v)$ of the coadjoint orbit 
$\mu(G(m))$ at $v =\mu(m)$
is the product $G_m \cdot S^1_m$, where $G_m$ is the stabilizer of the orbit
 $S^1_m$ is a subgroup in $G$, generating
the flows $\mu^{-1}\{ \mu(m)\}$ (Lemma 2.2). 
More precisely, since $Z(v)$ is connected and $dim \, S^1_m =1$, $Z(v)$
is the ``almost'' direct product of the connected component $G_m^0$ of $G_m$ with $S^1_m$. Here ``almost'' means that on the level of Lie algebras the
product is direct, and hence $G^0_m$ intersects with $S^1_m$ at a
finite group $\Z_p^0$.




\medskip


By Proposition 2.3  there  are  two     singular orbits $G/G_{min}$ and $G/G_{max}$ in  $M$.  
Fix a  geodesic segment  $\delta$ on  $M$ (we refer  the reader  to \cite{A-A}  and \cite{AA1993} for  discussion of  the notion of geodesic segment).  Denote by\\  
- $Z(v)$ the   stabilizer of $\mu(G(m))$, $m \in  \delta \cap   (M \setminus  (G/G_{min} \cup G/G_{max}))$,\\
- $Z_{min}$ the stabilizer of   $\mu(G(m'))$, $m ' \in \delta \cap  G/G_{min}$,\\
- $Z_{max}$ the stabilizer of   $\mu(G(m''))$, $m'' \in \delta \cap  G/G_{max}$.

The following  Lemma is straightforward.

{\bf Lemma  2.6}.  {\it There are only four possible cases:}

I) $Z(v) \cong Z_{min} \cong Z_{max}$

II) $Z(v) \cong Z_{min} \subset Z_{max}$

III) $Z(v) \cong Z_{max} \subset Z_{min}$.

IV) $Z_{max} \supset Z(v) \subset Z_{min}$


\medskip

Now we shall describe $M$ according to  four  cases 
in Lemma 2.6.

 Case (I):  all symplectic quotients $G(m) /S^1$ are $G$-diffeomorphic. In this case by dimension reason and the fact
that $G/Z(v)$ is simply-connected, we see immediately
that a singular orbit $G(m')$ is $G$-diffeomorphic
to its image $\mu (G(m')) =G/Z(v)$.   

In this case     Proposition 2.3 in the previous  version  of this  paper holds. Namely we have

{\bf  Lemma  3.} (\cite[Proposition 2.3]{Le1998}) {\it  There is  a Hamiltonian  $S^1$-action on $M$ such that
 $G(\mu(m))$  is a symplectic quotient of $M$ under  this $S^1$-action.}

{\it   Proof}. (see also Lemma 2.4 in \cite{Le1998})   Let $H_G$  be  a (unique up to constant) $G$-invariant  function on $M$ which satisifes the following condition
$$2\pi||grad\,H_G|| = L( \mu^{-1} \{ \mu (m)\} ), \eqno(2.4)$$
where  $(,)$ denotes  the length of the  $G$-invariant metric. It is easy to verify that $H_G$ generates  the required   Hamiltonian  action. \footnote{In the case $\dim \mu^{-1}(m) \le 1$ for all $m \in M$  we   prove  that $\pi_1 (M) = 0$,   using (E1).  Hence the symplectic  vector field   generating  the $S^1$-action on $\mu^{-1} (m)$ is a Hamiltonian vector field.}\QED

To
specify the $G$-diffeomorphism type of $M$ it
is useful to use the notion of segment \cite{A-A}.  In our
 case we just consider  
 the gradient flow of the function $H_G$ on $M$.
 After a completion and a reparametrization
we get a geodesic segment $[s(t)], \; t\in [0,1]$, in
$M$ such that the stabilizer of all the interior
 point $ s(t), t\in (0,1)$,  coincide with, say, $G_m$. (We observe that both $[s(t)]$ and the geodesic through $m$ with the
initial vector $grad\, H_G (m)$ are characterized by the condition
that every point in them is a fixed point of $G_m$). Denote by $G_0$ and
 $G_1$ the stabilizers  at singular points $s(0)$ and
 $s(1)$. Looking at the image of
the gradient flow of $grad\, H_G$ under the moment
map $\mu$ we conclude that $G_0 =G_1$.


\medskip

{\bf Proposition 2.7.} {\it In the case (I) 
 $M$ is $G$-diffeomorphic to  $G \times _{G_0} S^2$,
where $G_0 = (G_m^0 \times S^1_m)/\Z^0_p$ is
 the almost direct product of $G_m^0$ and $ S^1_m$, and the left action of $G_0$ on $S^2$  is
obtained via the composition of the projections $G_0 \to 
S^1_m/\Z^0_p$ with a Hamiltonian action of $S^1_m/\Z^0_p$ on $S^2$. }

\medskip

{\it Proof.} First we identify
the singular orbits in $M$ and in $G\times_{G_0} S^2$. The segment $[s(t)]$ extends this diffeomorphism to 
a diffeomorphism between $M$ and $ G_{\times G_0} S^2$.
Since $H_G$ is $G$-invariant it follows that
this diffeomorphism is $G$-diffeomorphism.
 \QED


Now let us compute the cohomology ring $H^* (M, \R)$ ( for $M$ in the case I). 
Once we fix a Weyl chamber we get a canonical $G$-invariant projection $\Pi_{\mu}$:
$\mu(M) \to \mu (G(m_0))$, where $G(m_0)=G/G_0$ is a singular orbit
in $M$. Let $j:= \Pi_{\mu} \circ \mu$ denote the
projection $M \to B:=\mu(G(m_0))= G/Z(v) \cong G(m_0)$.
Geometrically $j(x)= j( \mu^{-1} ( \mu(x)))$ is the limit of the flow generated by $grad\, H_G$  passing through
$x$.   Note
that $G(m_0)$ is the image of a section $s: B \to M$
of our $S^2$-bundle, and in what follows we shall
identify the base $B$ with ist section $G(m_0)$.  Let $f$ denote
the Poincare dual to the homology class $[G(m_0)]\in H_*(M,\R)$. 

Let $x_0 \in
H^2(\mu(G(m_0)), \R)$ be
the image of
the Chern class of the $S^1$-bundle $G(m) \to G(m_0)$, where
$G(m)$ is a regular orbit $G/G_m$ (or in other words, $x_0$ is
 the Chern class of the normal bundle over $G(m_0)$ with the
induced (almost) complex structure). 

Let $\{ x_i, R_1\}$ denote
the set of generators and their relations in cohomology ring $H^* (\mu(G(m_0)), \R)$ (see \cite{Bo}, correspondingly Proposition A.4 in  Appendix  A).

\medskip

{\bf Proposition 2.8}. {\it We have the following
isomorphism of additive groups
$$ H^*( G\times _{G_0} S^2, \R) = H^* (G/G_0, \R) \otimes H^*(S^2, \R).\eqno (2.5)$$

The only non-trivial relation in the algebra
$H^*(M, \R)$ are $R1$, $R2$, with
$$  f (f-j^*(x_0)) =0. \eqno (R2)$$}


{\it Proof}. The statement (2.5) on the additive structure of
$H^*(M, \R)$ follows from
the triviality of the cohomology spectral sequence of
our $S^2$-bundle. Clearly $(R1)$ remains the relation between 
the generators $\{j^*(x_i)\}$ in $H^*(M, \R)$. To show that the relation $(R2)$ holds we have two arguments.
One is in the proof of Lemma 2.12 and the other is here.
Using the intersection formula for $x_0$
we notice that the restriction
of $(f-j^*(x_0))$ to $G(m_0)$ is trivial. Thus to get
the relation ($R2$) it is enough to verify  that the value
of the LHS of $(R2)$ on the 
cycles in $M$ of the forms $j^{-1} ([C])$ is always zero, where $[C] \in H_2(B, \Z)$. Denote by $PD_M (.)$ the Poincare dual in 
$M$.
From  the identity
 $$ (PD _M( [ G(m_0)]) )^2 = PD_M ([G(m_0) \cap G(m_0
)])
=PD_M [ PD_B( x_0)] $$
we get $f^2 = PD_M [ PD_B (x_0)]$. Now it follows that 
$$f^2(j^{-1} ([C])) = f([C]).\eqno (R2.a)$$ 
On the other hand, since
the restriction of the 2-form representing $j^*(x_0)$ to
the fiber $S^2$ is vanished, we can apply the Fubini formula 
to the integration of a differential form 
representing  the class $f \cdot j^*(x_0)$,
(we can assume that $[C]$ is   represented by a pseudo manifold). In the result we
get  that
$$f\cdot j^* (x_0) (j^{-1} ( [ C]))= x_0 ([C]) )= f ([C]). \eqno (R2. b)$$ 
Thus (R2) is a relation in $H^*(M,\R)$. Finally the statement that (R2) is the only ``new'' relation in $H^*(M,\R)$ follows from the triviality of
our spectral sequence. \QED

{\bf Remark 2.9}.  If we take the other singular orbit $G(m_1) = G/G_1$
then  the Chern class of the $S^1$-bundle : $G(m) \to G(m_1)$ is
$-x_0$ (after an obvious identification
$G(m_0)$ with $G(m_1)$ since $G(m_1)$ can be considered
as another section (at infinity)
of our $S^2$-bundle). It is also easy
to see that the restriction of
$f$ on $G(m_1)$ is zero since $G(m_0)$ has no common point
with $G(m_1)$.


\medskip

{\bf Proposition 2.10}. {\it Let $M ^{2n}$ be in the case I of Lemma 2.6
and let us keep the notation in Proposition 2.8 for $M$. Then $M^{2n}$ admits a $G$-invariant symplectic form  $\om$ in a class $[\om] \in H^2 (M^{2n}, \R)$ if
and only if $[\om] =  j^*(x) + \alpha \cdot  f$ with $\alpha >0$, 
 and $(x+ t\cdot\alpha\cdot x_0)^{n-1} >0$ for all $t \in [0,1]$.  In particular $M^{2n}$ always admits
a $G$-invariant symplectic structure such that the action of $G$ on
$M$ is Hamiltonian. }

\medskip 

{\it Proof}.  Let $[\om] = j^*(x) + \alpha \cdot f$
with $ x\in H^2(G/G_0, \R)$.    The condition
that $\alpha>0$ follows from the fact that the
restriction of $\om$ to each fiber $S^2$ is positive. (Here
 we assume that the orientation
of $M$ agrees with that of $G(m)$ and the frame $(grad\, H_G, sgrad H_G)$.
The last frame is a frame of tangent space to the fiber $S^2$).
Thus the ``only if'' statement now follows trivially from the Duistermaat-Heckman Theorem.

Now let us
assume that the class $[\om]$ satisfies the condition
in Proposition 2.10.  Clearly
all these cohomology classes $ (x + t \cdot \alpha \cdot x_0)$, $ 0\le t \le 1$ are realized by
$G$-invariant symplectic forms by our condition (see also Remark A.5).
We fix a 1-parameter family of $G$-invariant metrics on $G/G_0$ which are also
compatible with these symplectic forms. According to Remark 2.9 (ii)
we can construct a  $G$-invariant metric on $M$ which compatible
with this family of $G$-invariant metrics on $G/G_0$.
Lifting to    $M$ we can define the restriction $\bar{\om}$ of $\om$ to each orbit $G(m)$. We normalize the $G$-invariant metric on
$M$ in the  direction $grad\, H_G$ orthogonal to the
orbit $G(m)$ such that the following condition (2.6) holds
$$grad\, H_0 (\bar{\om})(m) = - L (\mu^{-1}\{ \mu (m)\})j^* x_0, \eqno (2.6)$$
where $grad\, H_0 := grad\, H_G /||grad\, H_G||$ (we can 
normalize this metric by multiplying
 the length of $grad\, H_G$ with a positive function, because $\alpha >0$).
By the construction $\bar{\om}$ is a $G$-invariant 2-form on
$M$ whose rank is $(n-1)$. Denote by $\alpha\hat{f}_G$ the $G$-invariant 2-form on
$M$ whose restriction to each fiber $S^2$  is compatible with the restriction of the $G$-invariant metric 
to $S^2$.  We
put $\om = \bar{\om} + \alpha \hat{f}_G$. By the construction
$\om$ is a $G$-invariant 2-form of maximal rank on $M$. We claim that
$\om$ is a symplectic form realizing the class $j^*(x) + \alpha\cdot f$.
To verify the closedness of $\om$ it suffices to establish the
following identities
$$d\om (sgrad\, H_0, grad\, H_0, V_1) =0,\eqno (2.7)$$
$$d\om ( sgrad\, H_0, V_1,V_2) =0, \eqno(2.8)$$
$$d\om (grad\, H_0, V_1, V_2) =0, \eqno (2.9)$$
$$d\om (V_1, V_2, V_3) =0,\eqno (2.10)$$
for all $V_i$ in the normal bundle to the fiber $S^2$ and here  $sgrad\, H_0$ denotes the unite vector in $ker\, \bar{\om}_{| G(m)}$, whose orientation
agrees with that of the fiber $S^1$.
Using the formula
$$ 3d\om (X,Y,Z) = X (\om (Y,Z)) + Y(\om (Z,X)) + Z(\om (X,Y))$$
$$-\om ([X,Y], Z) - \om ([Y,Z], X) - \om ([Z,X], Y) \eqno (2.11)$$
we easily get that the LHS of (2.8) equals $d\bar{\om}_{| G(m)} =0$.

Applying (2.11) to (2.10) we also get that $ d\om (V_1, V_2, V_3) =
d\bar{\om} (V_1, V_2, V_3) + d\alpha \hat{f}_G (V_1, V_2, V_3) =0 +0 =0$.

To compute (2.7) we assume that $V_i$ is generated by
the action of a 1-parameter subgroup of $G$ (acting on $M$).  Taking
into account that $[sgrad\, H_0, grad\, H_0] \in ker\, \bar{\om}$ we get
$$ -3d\om (sgrad\, H_0, grad\, H_0, V) =\alpha\hat{f}_G ([grad \, H_0, V],
sgrad\, H_0)$$
$$ -\alpha \hat{f}_G ([sgrad\, H_0, V], grad\, H_0]).\eqno (2.12)$$
The RHS of (2.12) is zero since $\alpha\hat{f}_G$ is $G$-invariant. Hence (2.7) is zero.

To compute (2.9)  we also assume that $V_i$ is generated by
the action of a 1-parameter subgroup of $G$. Since
$H_G$ is $G$-invariant we get $[V_i, grad\,H_G] =0 = [V_i, grad\, H_0]$. Applying (2.11) to the LHS of (2.9) we get
$$3 d\om (grad\, H_0, V_1, V_2) = grad\, H_0 (\bar{\om} (V_1, V_2) ) -\alpha \hat{f}_G ( [V_1,V_2], grad\, H_0).\eqno (2.13)$$
   By the choice of $\alpha \hat{f}_G$ the second term in the RHS of
(2.13) equals\\
$-\la sgrad\, H_0, [V_1, V_2]\ra$.

Let us denote by $M^{reg}$ the set of regular points of the $G$-action on $M$. By the choice of $V_i$ and $\bar{\om}$ (see (2.6)) the first term in 
the LHS of (2.13) equals
${1\over 2\pi} d\theta (V_1, V_2)  \cdot L(\mu^{-1}\{ \mu (m)\})$
=$ - {1\over 4\pi}\theta([V_1,V_2])\cdot L(\mu^{-1}\{\mu(m)\} )$, where $\theta$ is the connection form on the
$S^1$-fibration $M^{reg}$.  In the presence of the (lifted) $S^1$-invariant
metric on $M^{reg}$ we can take $\theta ([V_1,V_2])$ as $4\pi \la sgrad\, H_0, [V_1,V_2]\ra /L(\mu^{-1}\{\mu (m)\})$. 

 It follows that
the LHS of (2.13) equals zero. This completes the proof of the closedness of $\om$.
Looking at the restriction of $\om$ to $G(m_1)$ and
$G(m_0)$ we conclude that $\om$ represents the class
$[j^*(x) + \alpha \cdot f]$.

\medskip

The statement on
the existence of a $G$-invariant symplectic
structure follows from the fact that $G/G_0$ always admits
 a class $x$ such that $x^{n-1} >0$. Since we can multiply $x$ with
a big positive constant $\lambda$, the class
$ (x+ t x_0)^{n-1}$ is also positive for all $ t\in [0,1]$ and
 we can apply the first statement here.

 The vanishing of the first
Betti-number of $M$ implies that the action
of $G$ is almost Hamiltonian and hence Hamiltonian because $G$ is compact.
\QED

\medskip

Cases (II) and (III) (in
Lemma 2.6). If we are
interested in the $G$-diffeomorphism type then these cases are equivalent.

\underline{Subcase (a):  $\dim \mu^{-1} ( m)  = 0$, if $m\in  G/G_{max}\cup G/G_{min}$}.
In this  subcase   the argument  in \cite{Le1998}, cases (II), (III)   is still valid.
W.l.o.g. it suffices to  consider the case (II): $Z(v) = Z_{min}$.
Clearly $G_{max} = Z_{max}$ and $G_{min} =Z_{min}$.  Note that  $G_{max}/G_{reg} = S^k$
by the slice theorem. On the other hand we have 
$Z(v) = G_{reg} \times S^1$.  Because $Z_{max}/Z(v)$ is
always of even dimension we  have $Z_{max}/Z(v) = \C P^{{k-1\over 2}}
=\C P^l$.

\medskip

{\bf Lemma 2.11}. {\it  In the case II, Subcase a,  we have the following
decompositions:  $ G_{max} = SU_{l+1} \times G_0$, $G_{reg}=
SU_l \times G_0$ and $Z_v = S(U_l\times U_1) \times G_0$, where the inclusion  $SU_l \to S(U_l \times U_1) \to SU_{l+1}$ is standard.}

\medskip

{\it Proof}. By checking the table A.3 (in the Appendix) of
 possible coadjoint orbit types we see that the
pair $(Z(v), Z_{max}\cong G_{max})$ in case (II a)
can be only:

Serie A. $Z_{max}= S(U_{l+1} \times \cdots \times U_{n_k})$.
Then $Z(v) = S(U_l \times S^1 \times \cdots  \times
U_{n_k})$ and $G_{reg} = S(U_l \times \cdots U_{n_k})$.

\medskip

Serie B, (D). $Z_{max} =U_{l+1} \times \cdots \times SO_{2n_k +(1)}$,  $Z(v)  = U_l \times U_1 \cdots  \times SO_{2n_k +(1)}$ and
$G_{reg}= U_n \times \cdots \times SO_{2h_k +(1)}$.

\medskip

Serie C.  Analogous to  B and D.

\medskip

Exceptional case: the same (see Table A.3 in Appendix). 

If $G$  is a product of compact Lie groups
then its coadjoint orbits are product of
coadjoint orbits of each factors. It it well-known that
every compact group Lie admits a finite covering which is
a product of compact simply-connected Lie group
whose algebra is simple. Thus to prove Lemma 2.11
in general case, it is not hard to see that it suffices to consider the above cases. 
\QED



\medskip

{\bf Proposition 2.12}.{ \it  Let M be in subcase (a) of case (II) (resp. of  case (III)). Then $M$ is $G$-diffeomorphic to   a $G$-invariant $\C P^{l+1}$-bundle over  $G/G_{max}$ (resp. $G/G_{min}$). There is a $G$-invariant  symplectic structure on $ M$ and
the action of $G$ is Hamiltonian with respect to this structure. }

\medskip

{\it Proof}. To prove the first
statement we consider the projection $M \to G/G_{max}:
\, x\mapsto \mu(x) \mapsto \Pi(\mu (x))$, where $ \Pi$ is a canonical projection
from $\mu (M)$ to the singular coadjoint orbit $G/G_{max}$. We
recall that this canonical projection can be chosen by using
the intersection
of $\mu(M)$ with a Weyl chamber (see \cite{Kir}).  By Lemma 2.11 the fiber of this projection
is the sum $D^{2(l+1)} \cup S^{2l+1}\times I \cup \C P^l$ and isomorphic
to $\C P^{l+1}$. Clearly this fiber consists of all trajectories of
 the flows $ grad\, H_G$ which end up at a point in the singular orbit
$G/G_{max}$. Hence the action of $G$ sends a fiber to a fiber.

It is also easy to describe the cohomology algebra
 of  $M$ by the method  in Proposition 2.8.
Namely we denote by $f$ the Poincare dual to the singular
 orbit $G/G_{min}$ of codimension 2 in $M$.
Since the singular orbit $G/G_{min}$ intersects the
fiber $\C P^{l+1}$ at a hyperplane $\C P^l$, the restriction
of $f$ on the fiber $\C P^{l+1}$ is the generator
of the cohomology group $H^2(\C P^n, \R)$.
Henceforth the ring $H^*(M, \R)$ is generated by $\{f, x_i\}$,
where $x_i$ are the pull-back of the generators
of the ring $H^*(G/G_{max} , \R)$ (compare (2.5)). Let $(R1)$ denote
the relation between $x_i$ in $H^*(G/G_{max}, \R)$, and
let $P_{min}$
 denote the Poincare dual to the singular orbit $G/G_{min} \subset M$. Put
$(R2)= f \cdot P_{min}$. It is easy to see (using the fact
 that two singular orbits have no common points and the associativity of the
cap action) that
(R1) and (R2) are the only relation
in $H^*(M,\R)$. (Now apply to the case in Proposition 2.8 we
observe that $P_{min} = f-x_0$). 

To
show the existence of
a $G$-invariant symplectic structure on $M$ we
use the lifting construction of
a family of invariant symplectic structures on $G/G_{max}$  as in the proof of Proposition 2.10. Here the main observation is the following.

{\bf Lemma 2.13}. {\it Let $G(m)$ be a principal orbit and $p_H$
denotes the projection from $M \setminus ( G/G_{max}) \to G/G_{min}$ which
 is defined
 by the gradient flow of $H_G$. Then
the characteristic leaf $ \mu ^{-1} \{ \mu (m)\}$  coincides
with $p_H ^{-1} (m) \cap G(m)$.}

\medskip

{\it Proof}. The projection of
the gradient flow of $H_G$ is
 also a gradient flow of a  $G$-invariant function
$H$ on $\mu (M)$. The slice theorem tells us that 
along the gradient flow of $H$ all the stabilizer groups  coincide.
Hence follows statement.\QED

Let  $[\om] = x + \alpha \cdot f$ be
an element in $H^2(M,\R)$. Clearly a necessary condition for
the existence of a symplectic form $\om$ in the class $[\om]$
is that $x^l >0$, $\alpha >0$ and for all $t \in [0,1)$ we have
that the restriction of the cohomology class $(j^*x+ t\cdot \alpha \cdot f)$
to the big orbit $G/G_{min}$ is also symplectic. (That follows from the Duistermaat-Heckman theorem or Kirwan's theorem).
 Here the restriction
of 
$f$ to the big orbit $G/G_{min}$ is the first Chern class of the $S^1$-fibration
$ G(m) \stackrel{ p_H}{\to} G/G_{min}$. Now let the class $[\om] \in H^2(M,\R)$ satisfy the above condition. Lifting the family of
symplectic forms on the quotient $(M \setminus (G/G_{max}))/S^1$
we get a symplectic form on $M\setminus (G/G_{max})$ ( see the proof of Proposition 2.10). By the construction
the lifted form extends continuously and non-degenerately
on the whole $M$
such that its restriction to the small orbit equals $j^*(x)$.
The closedness is also automatically valid. Considering the restriction
of the lifted form  to the two singular orbits
yields that our form realizes the cohomology class $ j^*(x) + \alpha \cdot f$.

To show the existence of a $G$-invariant symplectic 
structure on $M$ we use the fact that $G_{max}/G_{min}= \C P^l$. Under this
 condition we can find a $G$-invariant 2-form $\bar{x}$ in a class $x\in H^2(G/G_{max}, \R)$ such that $\bar{x}$ is a $G$-invariant
symplectic form and 
$j^*(\bar{x}) + t\bar{f}$ is a $G$-invariant symplectic form realizing
the cohomology class $ j^*(x) + t\cdot f$ for $t \in (0,1]$. (Here we construct a $G$-invariant 2-form on $G/G_{max}$ by $G$-invariant extension of a $G_{max}$-invariant
2-form  $\la \alpha, [X,Y]\ra $
in the $T_e(G/G_{max})$). 
This completes    our consideration of subcase (a) in cases (II) and (III).

\medskip

\underline{Cases (II) and  (III), subcase (b):  there is  $m \in G/G_{max}\cup G/G_{min}$ with $\dim \mu^{-1} (m) \ge 1$.}\\
W.l.o.g.  it suffices to consider   case (II). Since $Z_{min} = Z_v$,  using  the  relations $G_{min} \subset Z_{min}$ and  $G_{min} / G_m = S^l$,   there are only two possibilities
$$G_{min} = Z_{min}  = Z_v \LRA    G_{min} / G_m =  S^1,\eqno(E4)$$
$$G_{min}/G_{m} =  \Z_2. \eqno(E5)$$
The possibility (E5) cannot happen, since   in this case  $G_{min}$ and $G_m$   are defined  uniquely by the  $Z_{min} = Z_v$ and
the preimage $\mu^{-1} (x)$,   where $x \in G/G_{min}$    or $x \in G/G_m$ respectively,  are   connected  circles (the proof of Lemma 2.2 is also  valid for  exceptional  orbits).  Thus (E4) holds.
 In this case we have the following diagram of fibrations  and inclusions
$$\xymatrix{
&    G_{m}\ar@{^{(}->}@<3pt>[dl]\ar@{^{(}->}@<-3pt>[dr]& \\
G_{min} = Z_v  \ar[ur]^{S^1} \ar@{^{(}->}@<3pt>[dr]&  & G_{max}\ar[ul]^{S^n}\ar@{^{(}->}@<-3pt>[dl]\\
& Z_{max}\ar[ul]\ar[ur]^{\mu^{-1} (pt)} &
}$$

We  will call    a quintuple $(G, Z_v, Z_{max}, G_{max}, G_m)$    {\it admissbile}, if \\
- $G$ is  a  connected compact  group and $Z_v, Z_{max}, G_{max}, G_m$ are its  compact  subgroups,  moreover $Z_v \subset Z_{max}$  are   stabilizers  of  coadjoint orbits  in  $Lie\,  G$,\\
-  $Z_{max} \not = G_{max}$, and $Z_{v} = G_m \cdot   S^1 $,  $G_{max}/ G_m =  S^n$.

An admissible quintuple $(G, Z_v, Z_{max}, G_{max}, G_m)$  will be called {\it effective}, if  there are  no normal subgroups $G_1, G_2$  of $G$  and   a stabilizer  $H_1\subset G_1 $ of  a    coadjoint orbit
in $(Lie\, G_1)^*$   such that 
$$ G = G_1 \cdot G_2, \, Z_v =  H_1 \cdot  Z^2_v,\,   Z_{max} = H_1 \cdot  Z^2 _{max}, \, G_{max} = H_1 \cdot G^2_{max},  \, G_m = H_1 \cdot G^2_m$$
 and
$(G_2, Z^2_v, Z^2_{max},  G^2_{max}, G_m)$   is an admissible  quintuple.



Beginning  with the list of all possible     stabilizers  $Z_v \subset Z_{max}$  of  the coadjoint orbits of a compact Lie group $G$ (Table A.3 in Appendix A)  we pick up  from them  all possible
 triples $(G_m \subset  Z_v\subset G)$  and $(G_{max} \subset Z_{max}\subset G)$ such that  $G_m\cdot S^1 = Z_v$,   $G_{max} /G_m = S^n$ and $G_{max} \not  = Z_{max}$  with the help of   the table  of representation   of the sphere $S^n$  as an effective homogeneous space due to  Montgomery-Samelson-Borel and compiled by Alexseevsky-Alekssevsky  \cite[Table 1]{AA1993}, see also \cite{Borel}.  As a result we  compile   the following list of all    effective  admissible   quintuples.  

1) $(G= SO(2k+1), Z_{v} = SO(2k-1) \times S^1, Z_{max} = G= SO(2k+1),  G_{max} = SO(2k), G_m = SO(2k -1))$,  $k\ge 2$,\\
1a) $(G = SO(3), Z_v =  S^1,  Z_{max}  = SO(3), G_{max}  = SO(2),  G_m = \Z_p$),\\
2)$(G =  SO(2k+2), Z_{v} = SO(2k) \times S^1, Z_{max} =G = SO(2k+2),  G_{max} = SO(2k+1), G_m = SO(2k)$, $k \ge 1$,\\
3)$(G= SU(3), Z_v = U(1), Z_{max} = G= SU(2), G_{max} = U(1),  G_m = \Z_p)$.

In the  cases (1), (2), taking into account  the  above diagram,  we conclude  that 
for each $n\ge 3$  there is  a unique    effective  admissible    quintuple $(G= SO(n+1), Z_{v} = SO(n-1) \times S^1, Z_{max} = G= SO(n+1),  G_{max} = SO(n), G_m = SO(n -1))$.   It is not hard to see  that the corresponding 
 compact symplectic    manifold   that admits     cohomogeneity  1  Hamiltonian $SO(n+1)$-action   is    the  Grassmanians of   oriented 2-planes
$SO(n +2)/ (SO(n) \times SO(2)$   provided  with $SO(n+1)$-actions via the   standard  inclusion  $SO(n+1) \to SO(n+2)$, see e.g. \cite{Audin}.  By     the Alekseevsky-Alekseevsky  theorem \cite[Theorem 7.1]{AA1993}  the  above manifolds  are    the only ones  (up to $G$-diffeomorphism) that admit  a $G$-action of cohomogeneity 1 whose   orbit  types  are  listed  above.
The cases (1a)  and  (3)   correspond  to   4-dimensional symplectic manifolds   and they are well understood  \cite{I, Au}.  
The case (1a)  corresponds  to  the  action of $SO(3)$ on  $\C  P^2$ via  the embedding $SO(3)\to  SU(3)$.  In the case  (3)   the corresponding manifold is $\C P ^1 \times \C P^1$ with the diagonal action  of $SU(2)$.  Hence  we  obtain the following.

{\bf Proposition 4.} {\it  Let  $M$  be a        compact differentiable  $G$-manifold  of cohomogeneity 1      corresponding   to  one of the cases
listed above. Then   $M$ admits a  symplectic  form  which is $G$-invariant.}


This completes  our   consideration in cases I, II, III.

Now let us  consider case (IV).  

{\bf Lemma  5.} {\it  In  case (IV)   we have $\dim  \mu ^{-1} (m) \le 1$ for  all  $m \in  M$.}

{\it Proof.}  
Assume  the opposite,  w.l.o.g.  we can assume that  $\dim \mu^{-1} (m) \ge 1$, if $m \in G/G_{max}$.
Then   the   quintuple $(G, Z_v, Z_{max}, G_{max}, G_m)$  is  admissible.   Above,  we have classified in  Case (II), subcase b,  all  effective  admissible  quintuples.



In cases  (1),  (2),   the conditions $  G \supset Z_{min} \not =  Z_{v}$  and $Z_{min} \supset Z_v$ imply  that $Z_{min}  = Z_{max} = G$. Taking into account $G_{min}/G_m = S^n$ we conclude that  $G_{min} = G_{max}$.  Since $Z_{min}  = Z_{max} = G$, the
 singular orbits are   Lagrangian  spheres $G/G_{min}$ and $G/G_{max}$.  Using (E1),  we conclude  that  there exists  a nonzero constant $c$ such that the  symplectic  form $\om $ on $G(m) \times   (0,1) \subset M$ has the following  form
$$ \om (t, y) =  c\cdot ( dt \wedge \alpha  + t d\alpha) \eqno (E6)$$
where  $ t \in (0,1)$ and   $\alpha$ is the    canonical  connection 1-form of the $S^1$-bundle $G(m) \to   G(\mu (m))$.

Let $m_0$ and  $m_1$  be two  points on the      singular orbits  corresponding to $G_{min} $  and $G_{max}$.
By Weinstein theorem  the  neighborhoods  $U(G (m_0))$  of $G(m_0)$ and  $U(G(m_1))$ of  $G(m_1)$ are symplectomorphic.    Now (E6) implies that
 $G(m) \times  (0,1)$  cannot glue  with both   $U (G(m_0))$ and  $U(G(m_1))$, preserving the symplectic form $\om$. (By cohomological consideration this is possible only of $\dim M \le 4$. This  dimension  has been    considered in \cite{Au}.)   This  completes the proof
of  Lemma 5.\QED


The same argument as in case (II),  subcase (a), \footnote{since  Lemma 5 holds,    the argument in the  original  version of this  paper  is  valid}  shows that $G_{max} \cong Z_{max}$, $G_{min} \cong Z_{min}$
and $Z_{max}/Z(v)=\C P^l$, $Z_{min}/Z(v)= \C P^k$. 

\medskip

{\bf Proposition 2.14}. {\it Suppose that $M$ is in case IV. Then
$M$ is $G$-diffeomorphic to a $G$-invariant $\C P^k$-bundle over a coadjoint
orbit of $G$ or to the symplectic blow-down of such a $G$-bundle
along  the two singular (simplectic) orbits of $G$.}
\medskip

\Proof. We consider 3 possible subcases: (IVa), (IVb), (IVc).

\medskip

(IVa) \underline{If $ l\ge 2$ and $k \ge 2$}, then $G_{max} = S(U_{l+1}\times U_k\times U_1)\times G_0$,
 $G_{min}=S(U_l\times U_1\times U_{k+1})\times  G_0$, $G_{reg}=
S(U_l \times U_k \times S^1) \times G_0$, and $Z(v)= S(U_l\times U_1 \times U_k \times U_1) \times G_0$. Here the inclusion
$U_l \to U_{l+1}$ and $U_k \to U_{k+1}$ is 
canonical. Let $\Oo: = G/ (S(U_{l+1} \times U_{k+1})\times G_0 )$ be a coadjoint orbit of $G$. Let $\Pi_{min}$ denote the natural $G$-equivariant
projection from $G/G_{min} \to \Oo$. In the same way we
define the projection $\Pi _{max}$. We observe that
if the two points $m_{max} \in G/G_{max}$ and $m_{min} \in G/G_{min}$
are in the same gradient flow of the $G$-invariant function $H_G$
then their image under $\Pi_{max}$ and $\Pi_{min}$ coincide.
Hence the projection $\Pi _{min} $ and $\Pi_{max}$ can be extended
to a  projection $\Pi: M \to \Oo$. Clearly
the fiber is invariant under the $G$-action.
The group $S(U_{l+1}\times U_{k+1})$ acts on  the fiber of projection $\Pi$
from $M$ to $\Oo$ with  three orbit types: the singular ones are $ \C P^l$ and $\C P^k$ and the regular orbit is $ S(U_{l+1}\times U_{k+1}))/S(U_l \times U_k \times S^1)$.
Thus the fiber is diffeomorphic to $\C P^{l+k+1}$.

The simplest example of this case is $\C P^{l+k+1}$ with the standard action by $S(U_{l+1} \times U_{k+1})\subset SU_{k+l+2}$.

\medskip

(IVb) \underline {If $k$ = 1, $l \ge 2$},  then except the
above decomposition for $G_{max}, G_{min}$, $G_{reg}$ and $Z(v)$ there is only 
the following  possible subcase: $Z(v) = S(U_1 \times U_1 \times U_l ) \times G_0$, $G_{max} = S(U_2 \times U_l) \times G_0$, $ G_{min} = S(U_1 \times U_{l+1}) \times G_0$,
and $G_{reg}= SU_l \times S^1 \times G_0$. Let $S^1_m$ be
the subgroup of $Z(v)$ generated by
the vector orthogonal to $Lie\, G_{reg}$ in $Lie\, Z(v)$. Denote
by $\tilde{M}$ the suspension of $G/G_{reg}$. Clearly $\tilde{M}$ is
diffeomorphic to $G \times_{Z(v)}S^2$, where 
$Z(v)$ acts on $S^2$ via the projection to $S^1_m$. According
to Proposition 2.10 $\tilde{M}$ can be provided with a
$G$-invariant symplectic form such that the reduced symplectic form at
$G/Z(v)$ (considered at the ``mean point'' in $\tilde{M}$) is
the same as that reduced from $M$. We claim that
$M$ is a symplectic blow down of $\tilde{M}$ along the two
singular orbits $G/Z(v)_{max}$ and $G/Z(v)_{min}$.
To see this we cut a $G$-invariant neighborhood of
two $G$-singular orbits in $M$ (resp. $\tilde{M}$). By
the very construction of $\tilde{M}$ these new symplectic manifolds
are symplectomorphic. Hence follows the statement.

Now we shall show the existence of such a $G$-symplectic
manifold. Denote by $k$ the Cartan subalgebra of $g$. By Kirwan's convexity theorem there are elements $v, \alpha \in k$ such that
$Z(v)= S(U_1\times U_1\times U_l) \times G_0$, $Z(v + \alpha)= G_{max}$,
$Z(v-\alpha) = G_{min}$. Duistermaat-Heckman tells us
that the Chern class of the $S^1_m$-bundle $G/G_{reg} \to G/Z(v)$
is proportional to $\alpha$. Hence the Lie subalgebra $Lie \, G_{reg}$ is orthogonal to $\alpha$ in $Lie\, Z(v)$. We shall show that
there are such elements $\alpha$ and $v$ satisfying the above condition.

 Without lost of generality we assume that
 $G_0=1$. Thus $G= SU_{l+2}$. Write $v=(x_1,x_2, x_3,\cdot _{l \, times}\cdot ,
 x_3)$ with $\sum x_i =0$ and $x_1\not = x_2$. Thus
 the equation for $\alpha = (\alpha_1, \alpha_2, \alpha_3, \cdots, \alpha_3)$
 is $\alpha_1 + \alpha_2 + l \alpha_3=0$, 
 $ x_2 + \alpha_2 = x_3 + \alpha_3$ (and  is not zero),
 $ x_1 - \alpha _1 = x_2 - \alpha_2$ (and is not zero).
 The solution to these equations is 
 $ (l+2) \alpha _1 = l(x_1- x_2)$,
$\alpha_2 =\alpha_1 - x_1 + x_2 = (l-1)x_1 - (2l-1) x_2$,
 $\alpha_3 = \alpha_2 + x_2 -x_3 = (l-1) (x_1 - 2 x_2) -x_3$.
 The only thing need to check is the fact that $Z_{max}/G_{reg} = S^{2l-1}$,
 $Z_{min}/G_{reg}= S^{2k-1}$, where $G_{reg}$ is the subgroup
generated by the subalgebra orthogonal to the vector $\alpha$. We can do
it by finding an orthogonal representation of $G_{min}$ (resp. $G_{max}$)
on $\C ^2$ (resp. $\C ^l$) such that it acts on $S^3$ (resp. $S^{2l-1}$)
transitively with $G_{reg}$ as an isotropy group (see also \cite{A-A} which includes
a corresponding Borel's table  of the groups transitively acting
on spheres).

With these data at hand it is easy to construct a $G$ invariant
symplectic structure on the  $G$-manifold ($G_{min}, G_{reg}, G_{max})$
by the same lifting construction as in the proof of Proposition 2.10.
Namely we chose the family of symplectic form on $G/Z( v+ t\alpha)$, $t\in [-1,1]$, as the Kirillov-Kostant-Sourriau form.

\medskip

(IVc) \underline{ If  $k=l=1$}, then except the decomposition  analogous in the   subcase (b) (and hence subcase (a)) there is only the following possible cases with $Lie\, G_{max}=Lie\, G_{min}= su_2\times Lie\, G_0$, $Lie\, Z_v = s(u_1 \times u_1) \times Lie\, G_0$. Using Kirwan's convexity theorem we conclude that this case never happens. \QED
\medskip

Clearly Theorem 1.1 follows from Lemma 2.6,   Propositions 2.10, 2.12, 2.14  and Proposition 4.
\medskip



\section { Small quantum cohomology
of some symplectic manifolds admitting a
Hamiltonian action with cohomogeneity 1
of $U_n$ .}

\medskip

 Small quantum cohomology\footnote{ for a definition
and a formal construction of  full quantum cohomology see \cite{K-M}} (or
more precisely the quantum  cup-product deformed at $
H^2(M,\C) \subset
H^* (M, \C)$) was first suggested by Witten in
context of quantum field theory and then has been
defined mathematically rigorous for semi-positive (weakly monotone) symplectic manifolds by Ruan-Tian \cite{R-T} (see
also \cite{M-S}) and
recently for all compact symplectic manifolds by \cite{F-O}.  This
quantum product structure is an important  deformation invariant
of symplectic manifolds (and recently M. Schwarz \cite{Sch} has
derived
a symplectic fixed points estimate in terms of
quantum cup-length). Nevertheless 
there are not so much examples of symplectic manifolds
whose quantum cohomology can be computed (see \cite{CF,  FGP, G-K, S-T, R-T, W}).  The
main difficulty in the computation of quantum cohomology is that if we
want to compute geometrically it
is not easy to ``see'' all the holomorphic spheres
realizing some given homology class in $H_2(M, \Z)$. (On the other hand, computational functorial relations for quantum cohomology are expected to be found).

 In this
section we consider only the case of $M$  being
a $\C P^k$-bundle over Grassmannian $Gr_k (N)$ of
 $k$-planes in $\C ^N$: $ M = U(N) 
\times _{(U(k)\times U(N-k), \phi)} \C P^k$,
where $\phi$ acts on $\C P^k$ through the composition
of the projection onto $U(k)$ with the embedding
$U(k) \to U(k+1)$ and the standard action of $U(k+1)$ on $\C P^k$ (``standard''
action means the projectivization of the standard linear action on $\C ^{k+1}$.)
It is easy to see that the action on $\C P^k$ of the restriction
of $\phi$ to $U(k)$ has two singular orbits: $\C P^{k-1}$ and
a point, and its regular orbits are the sphere $S^{2k-1}$. According to the previous
section we see that $M$ can be equipped with a $G$-invariant symplectic structure and a Hamiltonian
action of $G=U(N)$ with the generic
orbit of $G$-action on $M$ being isomorphic to $U(N)/(U(k-1) \times U(N-k))$
and its image under the moment map $\mu: M \to u(n)$ is symplectomorphic to
 the flag manifold $U(N)/(U(1) \times U(k-1) \times U(N-k))$. With respect to Lemma 2.6 we
see that $M$ belongs to the case (I) if and only
if $k=1$, in this
case $M$ is a toric manifold.  We can also
consider $M$ as the projectivization
of the rank $(k+1)$ complex vector bundle  over $Gr_k (N)$
which is the sum of the tautological $\C ^k$-bundle $T_0$
and the trivial bundle $\C$. A special case of such
$M$ is $\C P^2 \# \overline{\C P^2}$ whose quantum 
cohomology 
is computed in \cite[example 8.6]{R-T} (see also \cite{K-M}).

By Lemma 3.1 below $M$ admits a $G$-invariant monotone symplectic structure.
To compute the small quantum cohomology algebra of $M$ we use several tricks well-known before \cite{S-T,
R-T, W} (e.g. the use of Gromov's compactness theorem) and the positivity
of intersection of complex submanifold. (In our monotone case we
can also use the fact that the projection to the
base $Gr_k(N)$
of a holomorphic sphere in $M$ is also a holomorphic sphere
in $Gr_k(N)$ with area less or equal to the area of
the  original sphere). Thus we can solve this question in our cases positively. It
seems that by the same way  we can give a  recursive rigorous computation
of small quantum cohomology ring of full or partial flag varieties,  since
any $k$-flag manifold is a Grassmannian bundle over a $(k-1)$-flag manifold
(see also \cite{G-K, CF, FGP} for other approaches to this problem).

\medskip

Recall that \cite{Bo} the cohomology
algebra $H^*(Gr_k (N), \C)$ is isomorphic
to the factor-algebra of
the algebra $\C [ x_1, \cdots , x_k] \otimes \C [ y_{1}, \cdots, y_{N-k}]$ over the ideal
generated by $S^+_{U(N)} (x_1, \cdots , y_{N-k})$
(see also Proposition A.4 in Appendix A). Geometrically
$x_i$ is $i$-th  Chern class of 
the dual bundle of the tautological $\C ^k$-vector bundle over $Gr_k (N)$,
and  $y_i$ is $i$-th Chern class of the dual bundle
of the other complementary $\C ^{N-k}$-vector bundle
over $Gr_k (N)$. Another description of $H^*(Gr_k(N), \R)$
 uses Schubert cells which form an additive basis, the Schubert
classes, in $H^*(Gr_k(N), \R)$ (see e.g. \cite{FGP} and
the references therein for the relation
between two approaches). Summarizing we have (see e.g \cite{S-T,M-S})

$$ H^* (Gr_k (N) , \C) = {\C [ x_1, \cdots , x_k]\over
\la y_{N-k+1}, \cdots , y_N\ra} $$

where $y_{N-k+j}: = -\sum_{i=0}^{N-k+j} x_i y_{N-k+j -i}$ (are defined inductively).
The first Chern class of $T_*Gr_k (N)$ is
$N x_1$. 

The quantum cohomology of $Gr_k (N)$ was
 computed in \cite{S-T} and \cite{W}. Now let us
compute the quantum cohomology algebra  $QH^*(M, \C)$. Denote
by $f$ the Poincare dual to the big singular
orbit $U(N)/ ( U(1) \times U(k-1) \times U(N-k))$ in $M$. 
Let $x_1, \cdots , x_k$ be the generators of
$H^*(Gr_k (N), \C)$ as above.
It is easy to see  
that the first Chern class of $T_*M$ is $(N-1)x_1 + (k+1) f$.
Then the minimal Chern number of
$T_*M$ is GCD $(N  -1, k+1)$ 
(because the $H_2(M, \Z)$ is generated by $H_2(Gr_k(N))$
 and $H_2(\C P^k)$).

\medskip
{\bf Lemma 3.1}. {\it (i) We have
$$ H^*(M, \C) = { \C [f, x_1, \cdots , x_k] \over
\la f(f^k - x_1 f^{k-1} + \cdots + (-1)^k x_k), y_{N-k+1}, \cdots , y_N\ra}.$$
(ii) $M$ admits a $G$-invariant monotone symplectic structure.} 

\medskip

(i) The formula is known in more general context \cite[Chapter 4, \S 20]{B-T}, \cite[Chapter 4, \S 6]{G-H}. But in our simple case we shall
supply here a simple proof. To derive Lemma 3.1 from the proof of Proposition 2.12 it suffices to 
show that
$$PD_M (Gr_k(N) )= f^k - x_1 f^{k-1} + \cdots +(-1)^kx_k\eqno (3.1)$$
To prove (3.1)
we denote $PD_M (Gr_k(N))$ by a polynomial
$P_k(f, x_1, \dots , x_k)$. By considering the  restriction of $PD_M(Gr_k(N))$ to the small orbit
$Gr_k(N)$  we conclude that  the lowest term (free of $f$)
of $P_k$ is $(-1)^k x_k$.
To define
the other terms of $P_k$ we consider the restriction
of $PD_M (Gr_k(N))=P_k$ to the submanifold $\bar{M} \subset M$,
which is the $\C P^k$ bundle over $Gr_{k-1}(N-1)$. Let
$M'$ be a submanifold of $\bar{M}$ which is defined as $M$ but over
$Gr_{k-1}(N-1)$. Using the formula 
$$(P_k)_{|\bar{M}}= PD_{\bar{M}} (Gr_{k-1}(N-1)) = PD_{\bar{M}} (M') \cdot PD_{M'}{Gr_{k-1}(N-1)}  $$
and the fact that$ PD_{\bar{M}} (M') =f$,
we conclude (by using the induction step) that
$P_k$ equals  RHS of (3.1).

(ii) It is well-known that $Nx_1$ is a symplectic class in $H^2(Gr_k(N), \R)$. By checking the non-degeneracy of
the family of $U(N)$-invariant forms $(Nx_1 + t (k+1) f)$ 
at a point $T_e((U(N))/U(1) \times U(k-1) \times U(N-k)$
we conclude that the condition for the existence
of an invariant symplectic form in the proof of Proposition 2.12 holds. Hence
$M$ admits a $G$-invariant monotone symplectic structure.\QED
\medskip

According to a general principle for computing the small
quantum cohomology ring of a monotone symplectic
manifold $(M, \om)$ we need to compute
only the quantum relations (\cite{S-T, W}).
More precisely, let $g_i (z_1, \cdots, z_m)$ be
 polynomials generating the relations
ideal of the cohomology algebra $H^* (M, \C)$ generated
by $\{ z_i\}$. Then  $z_i$ are
also generators of the small quantum algebra $QH^*(M,\C)= H^*(M, \C) \otimes \Z [ q]$
with the new relations $\hat{g}_i (z_i) = q P_i (z_i, q)$. Here $q$ is the quantum variable, $\hat{g}_i$ is the polynomial defined by $g_i$ with
respect to the quantum product in $QH^*(M, \C)$.
Denote the quantum product by $\star$.    There are
several equivalent approachs to small quantum cohomology  but we use notations
(and formalism)  in \cite{M-S}. 
 
\medskip

{\bf Theorem 3.2}. {\it Let $M$ satisfy the
condition $ 2(k+1) = N-1$ and as before, let
$P_k$ denote the Poincare dual to $Gr_k(N)$. Then its
small quantum cohomology ring is isomorphic to}

$$ QH^* (M) = { \C [f, x_1, \cdots, x_k, q] \over
\la f\star P_k= q , y_{N-k+1}, \cdots ,y_{N-1}, y_N =(-1)^{k+1}q^2 f\ra  } $$

\medskip

{\it Proof}.  Recall that (see e.g. \cite{McD-S})
the moduli space $\Mm_A(M)$ of holomorphic spheres
realizing class $A\in H_2(M,\Z)$ gives a non-trivial
contribution the quantum product of $a \star b$, $a, b \in H^*(M, \C)$,
if there is an element $c \in H^*(M, \C)$ such that the Gromov-Witten-Invariant
$\Phi_A(PD(a), PD(b), PD(c)) \not=0$. In this case  we have
$$ deg\, (a) + deg\, (b) \le  dim\, M + 2c_1(A) \le \deg\, a + deg\, b + dim\, M, \eqno (3.2)$$
which is also called   a degree (dimension) condition.

Recall that in our case the minimal Chern number of $M$ is $(k+1)$.
Thus  from (3.2), Lemma 3.1 and the monotonicity condition we see  immediately that
if the moduli space $\Mm_A(M)$ has a non-trivial contribution
to the quantum relation then
$0< c_1(A) \le 2(k+1)$. Hence $A$ must be one of the five following
homology classes.

(C1) : the homology class $ [u]$ generating
the homology group $H_2(\C P^k, \Z) = \Z$ of
the fiber $\C P^k$; 

(C2) : class $ 2[u]$; 

 (C3) : class
$[v]$ which can be realized as a  holomorphic
sphere on one singular orbit $G(m_s)$ which is diffeomorphic to $Gr_k(N)$ (see also the previous section);

(C4) :  the (exceptional) class $[v]-[u]$,

(C5) : the (double exceptional) class $2([v]-[u])$.

 Note that
$[u] $ and $[v]$ are the generators of $H_2(M, \Z) = \Z \oplus \Z$.

Let us consider the moduli space of holomorphic spheres in
class $[u]$. It is easy to see
that with
respect to the standard integrable
complex structure $J$ on $M$ the $J$-holomorphic spheres realizing this
class $[u]$ are exactly the complex lines
of the fiber $\C P^k$. The simplest
way to see this is to look at the
projection of these holomorphic spheres on the
base $Gr_k(N)$. (It may be possible to see this
by using the curvature estimate in \cite{L}.  This curvature
estimate could be able to
show that the minimal sectional curvature
distribution in $M$ consists of 2-planes in
the tangent space of the fiber $\C P^k$. Using
the same curvature estimate we have characterized the
space of holomorphic spheres of minimal degree in 
complex Grassmannian and other complex
symmetric spaces \cite{L} as the space of
Helgason spheres.)   A simple computation
shows that the virtual dimension
of the moduli space $\Mm _u (\C P^1, M)$ of
$J$-holomorphic spheres realizing $[u]$
equals the real dimension of this space and equals $2(k+1) + 2k + 2 N(N-k)$.
We can also apply the regularity criterion $H^1(\C P^1, f^* (T_*M)) = H^1 (\C P^1, \bar{f}^*(T_*(\C P^k)) =0$.
Here $f$ is a $J$-holomorphic map $\C P ^1 \to M$ and
 $\bar{f}$ is its restriction  on the  fiber $\C P^k$. 

Now let us compute the contribution of the
moduli space $\Mm_{[u]}(M)$ to the quantum relations, i.e.   $A$ is in case (C1). First we note  that by dimension reason the quantum polynomial of degree
less than $(k+1)$ must  coincide with the usual polynomial
(in the ring $H^*(M,\C)$). Thus to compute the contribution
of $\Mm_{[u]}(M)$ to the first defining relation it suffices to compute
the following Gromov-Witten-Invariants with $1\le l \le k+1$
$$ \Phi_{[u]} (PD(f^l), PD (x_{k+1-l}), pt), \eqno (3.3a)$$
$$ \Phi_{[u]} (PD(f^{k}, PD(f), pt ). \eqno (3.3b)$$

We claim that  the Gromov-Witten-Invariant in (3.3a)
equals zero. 
We observe that $PD_M (x_{k+1-l}) = j^{-1}( PD_B (x_{k+1-l}))$,
where as in the previous section we denote by $j$ the
projection of $M$ to $Gr_k(N)$.
Hence, taking into account that $[u]$ is
 a ``fiber'' class we see immediately, by dimension reason,
that there is no 
holomorphic curve in class $[u]$ which intersects $j^{-1} (PD_B (x_{k+1-l}))$ and goes through a 
$PD (f^l)$.

We claim that  the G-W invariant in (3.3b) equals 1.
To prove this we fix a fiber $\C P^k$ which contains the given
point $pt$. We observe that  the singular orbit representing
$PD_M (f)$ intersects
with each fiber $\C P^k$ at a divisor $\C P^{k-1}$.
Finally we note that $PD_M (f^k)$ intersects with
the fixed $\C P ^k$ at one point because
$f^k ([\C P^{k-1}]) =1$. Since
there is exactly one complex
line through the given two points in $\C P^n$ (and
 this line
always intersects the divisor $\C P^{k-1} \subset \C P^k$)
we deduce that the G-W invariant in (3.3b) is 1.

Summarizing we get
$$ f\star_{[u]} P_k = q ,\eqno (3.3c)$$
(here the LHS of (3.3.3) denotes the quantum polynomial, deformed
by $[u]$).

\medskip

Next  we shall compute the contribution
of $\Mm_{[u]}$ to the ``old'' defining relation $y_j$, $ j = N-k+1,N$.
 First we shall show that
$$ \Phi _{[u]} ( PD_M (x_p), PD_M (y_{j-p}), PD_M [w])=0\eqno (3.3d) $$
for any $[w] \in H^*(M)$
with degree equal $dim\, M +  2(k+1)- 2j$. Using the formula $PD_M [j^* (y)] = j^{-1} PD_B [ y]$
for the Poincare dual of a  pull-back cohomology
class of the base of a fiber bundle we observe that
if (3.1) is not zero then $PD_M [w]) \cap PD_M (x_p) \cap 
PD_M (y_{j-p}) \not= \emptyset$. But it is impossible by
the dimension reason.

\medskip

Thus there remain possibly four other
non-trivial contributions,    associated with  cases (C2)-(C5),
to the quantum relations. The first one is related to the
Gromov-Witten invariants
$$ \Phi _{[2u]} ( PD_M (x_p), PD_M (y_{j-p}), 
PD_M (w)) , \eqno (3.4)$$

 the second to the Gromov-Witten invariants
$$\Phi_{[v]} (PD_M (x_p), PD_M (y_{j -p}), PD_M (w)),
 \eqno (3.5)$$

and the two  other Gromov-Witten invariants related to
the (exceptional) classes $[v] - [u]$ and $2([v]-[u])$.

Here, in the cases  (C2) and (C3),   the degree of $w$  in (3.4) and (3.5) must be
dim $M$ + $ 4(k+1) -2j$. 
 
To compute (3.4) we
use a generic almost complex structure $J_{reg}$ nearby
the integrable one. Thus
the image of $J_{reg}$-holomorphic
spheres in class $2[u]$ must in a (arbitrary)
small neighborhood of a complex line in the
fiber $\C P^k$, that is the
projection of
a $J_{reg}$-holomorphic sphere in class
$[u]$ must be in a ball of
 radius $\eps/2$. Now we can use the same
argument as before. Since $PD_M (x_p) \cap PD_M (y_{j-p}) \cap
PD_M (w) = \emptyset$ there exists a positive number $\eps$
such that the $\eps$-neighborhood of these cycles also
do not have a common point. Now looking
at the projection of these cycles on the base $Gr_k(N)$
we conclude that the
contribution (3.4) is zero.

\medskip

In order to  compute the contribution (3.5) we need
to know the moduli space of the holomorphic spheres in class
$[v]$ whose dimension is   dim $M$ + $4 (k+1)$ = dim $Gr_k (N)$
+ $ 6k +4$
= dim $Gr_k (N)$ + $2 N + 2(k-1)$. We pick up the standard integrable
complex structure. We claim that all these
holomorphic
spheres can be realized as holomorphic sections
of $\C P^k$-bundle over $\C P^1 _{[v]}$, where
$\C P^1 _{[v]}$ is a holomorphic sphere of
 minimal degree in $Gr_k (N)$. Indeed over
this $\C P^1$ the
bundle $\C P^k$ is
the projectivization
of the sum of $(k+1)$ holomorphic line
bundles with $k$ Chern numbers being 0 and
one number being $(-1)$. Thus for any holomorphic sphere $(S^2,f)$ which is
a holomorphic section of the $\C P^k$ bundle over $\C P^1$ we
have $H^1(S^2,f^* (T_*M)) = H^1 (S^2, f^*(T_*\C P^k)) =0$. To
show that these holomorphic sections exhaust all the  holomorphic
spheres in the class $[v]$ we look at
their projection on the base $Gr_k(N)$.
 
Now let us to compute (3.5) with $j = N-1$ or $j=N$
(by dimension condition (3.2) those are  the only cases  which may enter into the quantum relations).

If $j= N-1$ then   the contribution in (3.5)
must be 0 since we know that
on the base $B= Gr_k (N)$ there
is no holomorphic curve of minimal degree which
go through the cycle $PD_B (x_p)$ and
 $PD_B (y_{N-p-1})$ (by dimension reason). 

If $j = N$ then there are two possibilities for
$PD_M (w)$, namely they are $[u]$ and $[v]$ - the generators
of $H^2 (M, \Z)$.

 Let us consider the first case i.e.,  $PD_M (w)$ is
 a holomorphic sphere $u$ in the fiber $\C P^k$. The induction argument on
$Gr_k(N)$ (\cite{S-T,W}) shows that $p$ in (3.5) must
be $k$ and there is a unique  (up
to projection $j$) holomorphic sphere in class $[v]$
which intersects with $PD_M (x_k)$ and  $PD_M (y_{N-k})$
and satisfies  the following property: its image under the projection $j$ goes through the fixed point $j(u)\in Gr_k(N)$. 
Hence
we can reduce our computation of the corresponding contribution in (3.5)
to the related Gromov-Witten invariant
in the $\C P^k$-bundle over
$\C P^1_{[v]}$. Thus we get
$$ \Phi _{[v]} (PD_M (x_k), PD_M (y_{N
-k}), [u]) = (-1)^{k+1}. \eqno (3.5a)$$
 Now let us consider the second case i.e.,  $PD_M (w)$ is the class  $[v]$ realized by
a holomorphic section of the $\C P^k$-bundle over the $\C P^1$.
Clearly there is only one holomorphic section  passing through a given 
point in this bundle. Thus  we get
$$ \Phi _{[v]} (PD_M (x_k) , PD_M (y_{N-k}), [v])
 = (-1)^{k+1}. \eqno (3.5b)$$

\medskip

 Let us consider case (C4), i.e.  the moduli space of
holomorphic spheres in the class $[v]-[u]$. We have two arguments to show that there
is no $J$-holomorphic sphere in this class. The simplest argument was
suggested by Kaoru Ono. Namely considering the
intersection of a holomorphic sphere in this class with the big singular
orbit $U(N)/( U(1) \times U(k-1) \times U(N-k))$ yields that there
is no holomorphic sphere in this class.  The another (longer) argument
uses the area comparison. 
Clearly the area of such a holomorphic sphere equals the value $\om ([v] -[u])$. On the
other hand the projection to $Gr_k(N)$ of a holomorphic sphere in
this class has area $\om ([v]) > \om ([v]-[u])$. ( The projection decreases
the area  because of Duistermaat-Heckman theorem applied to our
monotone case). Thus
there is no $J$-holomorphic sphere in this class.
Since the class $[u] -[v]$ is indecomposable
in the Gromov sense it follows from
the Gromov compactness theorem that for
nearby generic almost complex structure $J'_{reg}$ there
is also no $J'_{reg}$-holomorphic sphere. Thus there  is no
quantum contribution of this class.

\medskip

 Finally we consider the quantum contribution in case  (C4) associated to the class $2([v]-[u])$. The space of $J$-holomorphic spheres in this class is empty by the
same reason as above (two arguments).  Finally by using the Gromov compactness
theorem we can show  the existence of
a regular almost complex structure $J_{reg}$ nearby $J$ such
that there is no $J_{reg}$-holomorphic sphere
in this class. (Because if bubbling happens,
they must be holomorphic spheres in class $[v]-[u]$, which is also impossible.) 

\medskip

Summarizing we get that
the only new quantum relations are those involving
(3.3c), (3.5a) and (3.5b).  Note that
$f$ is defined uniquely by the condition $f(u) = 1 = f(v)$.
This completes the proof of  Theorem 3.1.\QED

{\bf Remark 3.3.} Since the rank of $H_2(M)$ is
2 it is more convenient to take 2 quantum variables $q_1, q_2$.
In this case our computations  give a (slightly)
formal different  answer, namely $(R2)= q_1$
and $y_N = (-1) ^{k+1} (q_1^2f_1 + q_2^2 f_2)$.
Here $f_1$ and $f_2$ form a basis of $Hom\, (H_2 (M, \C), \C)= H^2 (M, \C)$ which is dual to the basis $([u],[v])\in
H_2(M,\C)$.

\medskip

{\bf Remark 3.4}. Let $M$ be a symplectic manifold as
in Theorem 3.2. 

(i) It follows immediately from
Theorem 3.1 and Schwarz's result \cite{Sch} that
the any exact symplectomorphism on
$M$ has at least $k+1$ fixed points.

(ii) It seems that after a little work we can apply
the result in \cite{H-V} to show
that the Weinstein conjecture also holds for
those $M$.

\bigskip

\section{ Compact symplectic manifolds
admitting symplectic action of cohomogeneity 2}

\medskip

A direct product of $(M_1, \om_1)$ and
$(M_2, \om _2)$ is
 a symplectic manifold which admits a symplectic
action of cohomogeneity 2 provided that
either both $(M_i, \om_i)$ admit symplectic
action of cohomogeneity 1 or $(M_1, \om_1)$ 
is a homogeneous symplectic manifold
and $(M_2, \om_2)$ has dimension 2. These examples
are extremally opposite in a sense that,
in the first case the normal bundle
of any regular orbit is isotropic, and
in the second case the normal bundle is symplectic.

\medskip

{\bf Proposition 4.1}. {\it Suppose that an action
of $G$ on $(M^{2n},\om)$ is Hamiltonian 
and $dim\, M/G =2$. Then either all the principal
orbits of $G$ are symplectic (simultaneously),
or all the principal
orbits of $G$ are coisotropic (simultaneously). 
In the first case a principal orbit is isomorphic to
a coadjoint orbit of $G$, in the last case a principal
orbit must be a $T^2$-bundle over
a coadjoint orbit of $G$.}

\medskip

{\it Proof}. Since the set $M^{reg}$ of
regular points in $M^{2n}$ is
 open and dense in $M^{2n}$, and the property of being symplectic is an open
condition, it suffices to show that there is an open, dense, $G$-invariant set 
$ M^{\circ} \subset M^{reg}$ such that all the
orbit $G(x) \subset M^{\circ}$ is symplectic (or coisotropic simultaneously).
We consider the moment map  $\mu: M^{2n} \to g^* =g$.
By  Sard's theorem  the set $S_{\mu}$ of points $x$ in $M^{2n}$, where
the dimension $d$ of $\mu^{-1} \{ \mu(x)\}$ is maximal, is open and dense in $M^{2n}$. Let $ M^{reg}_{\mu}$ be the  set in $M$ consists
of points $x$ such that $\mu(G(x))$ is a orbit of maximal dimension in
$\mu(M)$. Using Kirwan's theorem we see that $M^{reg}_{\mu}$ is an open and dense set in $M$.
We claim that we can take  $M^{\circ}$ as the intersection
of $S_{\mu}$ with $M^{reg}_{\mu}$ and the set of
regular points in $M^{2n}$.  Using the formula (2.3) we note that $d\le 2$.   Since the dimensions of $G(x)$ 
  and  of $\mu (G(x))$ are even  if $x\in M^{\circ}$,   we get that $d$ must be either 0 or 2. First we suppose  that
$d= 0$. Since $G$ is connected all the
other principal orbit $G(m')$ in $M$ also connected, and
since $\mu(G(m))$ is simply connected,
all the principal orbits in $M^{\circ}$ must be diffeomorphic
to $\mu (G(m))$ (and hence 
are symplectic). Clearly if
 orbit is symplectic then the restriction of $G$-action on it is also Hamiltonian, thus by Kirillov-Kostant-Sourriau theorem, it
must be isomorphic to a coadjoint orbit of $G$. Now let us
assume that the
 ``generic'' dimension $d$ of $\mu ^{-1} \{ \mu (m)\}$
is 2. Since the dimension
of $\mu(G(x))$ is a constant for  $x \in M^{\circ}$, we conclude that
either all $G(x)$, for $x\in M^{\circ}$ is  either symplectic  simultaneously
or isotropic simultaneously. In the last  case $\mu^{-1} \{ \mu(x)\} \subset G(x)$ and
 $ \mu(G(x)) = G(x) / \mu^{-1} \{ \mu(x)\}$.
 Arguing as in the proof
of  Proposition 2.1
 we see that  $\mu^{-1}\{ \mu (x)\}$ admits a
nowhere zero vector fields $sgrad\, \Ff _{v_1}$  and $sgrad \, \Ff _{v_2}$. Thus
it must be an isotropic torus. \QED

{\bf Remark 4.2}. {\it  (i). The quotient space
$\mu(M)/G$ is either a point or a convex 2-dimensional
polytope.

(ii) If the action
of $G$ is Hamiltonian and the principal orbit
is symplectic then the condition
that $\mu(M)/G$ is
a point is equivalent to the
fact that  $d$ (in the proof of Proposition 4.1) equals 2. In this
case $M$ is diffeomorphic to a bundle over a coadjoint orbit of $G$ whose
 fiber is a 2-dimensional surface.}

\medskip

 The first statement in Remark 4.2 follows from the proof of Proposition 4.1
and Kirwan's theorem on convexity of
moment map. The second statement follows by considering
the moment map.

\medskip

{\bf Proposition 4.3}. {\it  Suppose that the action
of $G$ is Hamiltonian, the number $d$ (in the proof of Proposition  4.1) is zero and the action of $G$ on $\mu(M)$ has only one orbit type. Then $M$ is
$G$-diffeomorphic to a  fiber bundle
over a 2-dimensional surface $\Sigma$, whose fiber is isomorphic to
a coadjoint orbit of $G$.}

\medskip

Indeed, by the dimension reason in
this case there is also only one orbit type of $G$-action on $M$. Note that such a bundle
always admits a $G$-invariant symplectic structure. 

\medskip

If the principal orbits of $G$ in $M$ are coisotropic then
$P=\mu(M)/G$ is  a 2-dimensional convex polytope.

\medskip

{\bf Proposition 4.4}. {\it If
the action of $G$ on
$M$ is Hamiltonian and
the principal orbit of $G$ is coisotropic
then $M$ is diffeomorphic to the bundle of
ruled surface over a coadjoint orbit
of G provided that the action
of $G$ on $\mu(M)$ has only one orbit type.}

\medskip

{\it Proof}. In this case $M$ admits a projection $\pi$ over
 a coadjoint orbit $\mu(G(m))$ with fiber $\pi^{-1}$ being
a symplectic 4-manifold. This symplectic 4-manifold admits
a $T^2$-Hamiltonian action. Hence it must be a rational or ruled surface (see \cite{Au}). \QED

\bigskip

{\bf Appendix  A.
 Homogeneous symplectic spaces of compact Lie groups.}

\medskip

First we recall
a theorem of Kirillov-Kostant-Sourriau (see e.g. \cite{Kir}).
\medskip

{\bf Theorem A.1.}  {\it A  symplectic
manifold admitting a Hamiltonian homogeneous action of a connected Lie group $G$ is isomorphic to
a covering of a coadjoint orbit of $G$.}

\medskip

If $G$ is a connected compact Lie group,{\it  using the homotopy exact sequences}, it  is not hard to see that
all its coadjoint orbits are simply-connected.
Thus in this case we have the following simple

\medskip

{\bf Corollary A.2.} {\it A symplectic
manifold admitting a Hamiltonian homogeneous action of a connected compact Lie group $G$ is a
coadjoint orbit of $G$}. 

\medskip

{\bf Table A.3.} We present here a list of all coadjoint orbits
of simple compact Lie groups. Recall that a coadjoint orbit
through $v\in g$ can be identified with the homogeneous space $G/Z(v)$ with $Z(v)$ being the centralizer
of $v$ in $G$. Element $v$ in
a Cartan algebra $Lie\, T^k \subset g$ is regular iff for all root $\alpha$ of
$g$ we have $\alpha (v) \not= 0$. In this case $Z(v)$ is the  maximal torus $T^k$
 of $G$. If $v$ is a singular element with
$\alpha_i (v) =0$ then $Lie \, Z(v)$ is a direct sum
of the subalgebra in
$g$ generated by the roots
$\alpha _i$ and $Lie\, T^k$. To identify the type
of this subalgebra $Lie\, Z(v)$  we observe that $Lie\, T^k$ is
 its Cartan subalgebra and the root system of $Lie\, Z(v)$ consists of
those roots $\alpha $ of $G$  such that $\alpha (v) =0$. Looking
at tables of roots of simple Lie algebras \cite{O-V} 
and their Dynkin schemes we
get easily the following list (which perhaps could be found somewhere else)

(A). If $G = SU_{n+1}$ then $Z(v) = S (U_{n_i} \times \cdots
\times U_{n_k})$, $\sum n_i = n+1$.

(B,C,D). If $G$ is in $B_n$, $Z_n$ or $D_n$ then $Z(v)$
is a direct product $U_{n_1} \times \cdots U_{n_k}\times G_p$
with $rk G_p + \sum n_i = rk G$, and $G_p$ and $G$ must be
from the same series $B$, $C$, $D$.

Analogously but more combinatorically complicated
are the types of $Z(v)$ in the exceptional series. Note that
all the listed below simple exceptional groups are simply connected.

($E_6$). Except the regular orbits with $Z(v) = T^6$ we also
have other possible singular orbits  with $Z(v) = S(U_{k_1}\times\cdots \times  U_{k_n})$ with $n \ge 2$, $\sum k_i = 7$ and
$ T^{k} \times Spin _{6-k}$ with
$ k= 1, 2$.

($E_7$). Analogously. Possible are also $Z(v) = T^1 \times SU_2 \times
Spin_{10}$ and $Z(v) = T^1 \times E_6$.

($E_8$). Analogously. ( Possible are also $T^1 \times E_7$ and
$T^1 \times SU_2 \times E_6$).

($F_4$). Singular orbits can have $Z(v)$ being $T^2\times SU_3$, $T^2\times SU_2 \times SU_2$ or $T^1 \times Spin_7$ and $ T^1 \times Sp _3$. 

($G_2$) Except the regular orbit $G_2/T^2$ there are
also singular orbit $G_2/SU_2 \times T^1$.

\medskip

To compute the cohomology
ring of $G/Z(v)$ we
use:

\medskip
{\bf Proposition A.4.} (\cite[Theorem 26.1]{Bo}). {\it
The cohomology algebra $H(G/Z(v), \R)$ is
a factor-algebra $S_{Z(v)}$ over the ideal generated
by $\rho^*_R (S^+_G)$ which equals the characteristic
subalgebra.

(ii) Let $s_1 -1, \cdots , s_l-1$ and
correspondingly, $r_1 -1, \cdots , r_l-1$ be
degree of the generators in $H^*(G)$ and $H^*(Z(v))$.
Then the Poincare polynomial of $G/Z(v)$ equals}
$$ {(1-t^{s_1}) \cdots (1-t^{s_l}) \over
(1-t^{r_1}) \cdots (1-t^{r_l}) }. $$

\medskip

Here $S_G$ is the algebra of $G$-invariant polynomials
in $g$ and $S_G^+$ is its subalgebra which is generated by monomials of
positive degree.

\medskip

{\bf Remark A 5}. All the $G$-invariant symplectic
form on $G/Z(v)$ are compatible
with the (obvious) $G$-invariant complex
structure. Thus all of them are deformation equivalent to a monotone symplectic form.



\medskip

{\bf Acknowledgement}. I   thank  Dmitri Alekseevsky, Michael
Grossberg,  Yuri Manin and Tien Zung Nguyen for stimulating and helpful discussions. My thanks are also due to Michel Audin
for pointing out an interesting paper by Iglesias [I]. Finally I would like to express my  sincere
gratitude to
Kaoru Ono for numerous critical comments and suggestions,
and I thank the referee  for useful comments.

{\it  I thank   Christopher    Woodward  for  pointing  out the gap  in my previous classification  and   for his interest  in the
corrected  version.}




\end{document}